\documentclass[11pt]{article}
\catcode`\@=11 \catcode`\@=12

\usepackage[utf8]{inputenc}
\usepackage{amsmath}
\usepackage{amssymb}
\usepackage{amsthm}
\usepackage{oldgerm}
\usepackage{multicol}
\usepackage{amsfonts}
\usepackage{amssymb}
\usepackage{latexsym}
\usepackage{amsthm}%custom qedsymbol
\usepackage{graphicx}
\usepackage{wasysym}
\usepackage{etex}%in order to make the package easybmat work
\usepackage{easybmat} %matrices 
\usepackage{blkarray} %matrices
\usepackage{microtype}
\usepackage{calc}%arithmetic operations allowed inside \setlength
\usepackage{cases}
\usepackage{multirow,bigdelim}%big braces

\usepackage[pdftex,colorlinks,urlcolor=blue,pdfstartview=FitH]{hyperref}

\newcommand{\Rm}{\mathbb{R}}

\newcommand{\mF}{\ensuremath{\mathcal{F}}}

\newcommand{\Nm}{\ensuremath{\mathbb{N}}}

\newcommand{\mA}{\ensuremath{\mathcal{A}}}

\newcommand{\vs}{\vspace{.2cm}}

\newtheorem{lem}{Lemma}
\newtheorem{thm}{Theorem}
\newtheorem{conj}[lem]{Conjecture}

\newtheorem{prop}{Proposition}

\newtheorem*{prop*}{Proposition}

\def\proof {\noindent{\sc{Proof. }}}
\def\qed {\mbox{}\hfill {\small \fbox{}} \\}
\def\lto{\longrightarrow}
\def\lmto{\longmapsto}

\def\leq{\leqslant}
\def\geq{\geqslant}

\newcommand{\deter}{{\rm det}\,}
\newcommand{\Span}{{\rm Span}\,}
\newcommand{\codim}{{\rm codim}\,}

\newcommand{\rank}{{\rm rank}\,}

\newcommand{\R}{\mathbb{R}}
\newcommand{\N}{\mathbb{N}}

\newcommand{\F}{\mathcal{F}}

\newcommand{\Rr}{\mathcal{R}}
\newcommand{\Cc}{\mathcal{C}}

\title{Some Remarks on  Thom's Transversality Theorem }
\author{Patrick  Bernard\footnote{membre de l'IUF} 
\hspace{.2cm}and Vito Mandorino}
\date{}
\begin{document}

\maketitle 

\vspace{1cm}
\begin{center}
-----
\end{center}
\begin{multicols}{2}

\begin{small}
\noindent
Patrick Bernard,
 Universit\'e  Paris-Dauphine,\\
CEREMADE, UMR CNRS 7534\\
Pl. du Mar\'echal de Lattre de Tassigny\\
75775 Paris Cedex 16,
France\\
\texttt{patrick.bernard@ceremade.dauphine.fr}\\

\noindent
Vito Mandorino,
 Universit\'e  Paris-Dauphine,\\
CEREMADE, UMR CNRS 7534\\
Pl. du Mar\'echal de Lattre de Tassigny\\
75775 Paris Cedex 16,
France\\
\texttt{mandorino@ceremade.dauphine.fr}\\

\end{small}

\end{multicols}
\vs
\thispagestyle{empty}
\begin{center}
-----
\end{center}

\textbf{Abstract. } 
We study Thom Transversality Theorem using a point of view, suggested by Gromov, 
which allows to avoid the use of Sard Theorem and gives finer informations
on the structure of the set of non-transverse maps.

\begin{center}
-----
\end{center}

\textbf{R\'esum\'e. }
On étudie le théorème de Transversalité de Thom 
en utilisant un point de vue, suggéré par Gromov, qui permet d'éviter 
l'usage du théorème de Sard et fournit une description plus fine de 
l'ensemble des applications non transverses.

\begin{center}
-----
\end{center}

\textbf{Sommario. }
Studiamo il teorema di trasversalit\`a di Thom utilizzando un punto di vista, suggerito da Gromov, che consente di evitare l'utilizzo del teorema di Sard e fornisce una descrizione pi\`u fine dell'insieme delle mappe non trasverse. 
\begin{center}
-----
\end{center}

MSC: 
\newpage
\section{Introduction}

It is well-known that ``most'' functions are Morse, which means 
that their critical points are non-degenerate.
Discussing this claim with some details
 will be an occasion to introduce and motivate the present work.
Let us fix some integer $r\geq 2$ and a dimension $d$.
Let $B^n$ be the open unit ball in $\Rm^n$, and $\bar B^n$ the closed unit ball.
We denote by $C^r(\bar B^n,\Rm)$ the space of functions which are $C^r$
on $B^n$, and whose differentials up to order $r$ extend by continuity to 
the closed ball $\bar B^n$. We endow it with the norm given as the sum of the supremums
of the differentials of order less than $r$. It is then a separable Banach space.
Let $F$ be an affine subspace of $C^r(\bar B^n,\Rm)$. In most cases $F$ will just be the whole space 
$C^r(\bar B^n,\Rm)$, but it is sometimes useful to consider finite-dimensional spaces $F$.
The map
\begin{align*}
e_1:B^n\times F& \lto \Rm^n\\
(x,f)&\lmto df(x)
\end{align*}
is $C^{r-1}$, and, when $F=C^r(\bar B^n,\Rm)$ it is a submersion, see \cite{AR}, Theorem 10.4.
Recall that a $C^1$ map is a submersion if and only if its differential
at each point is onto with a split kernel.
Then, it is locally equivalent (by left and right composition
by $C^1$ maps) to a  projection with split kernel.
We will always assume that $F$ is chosen such that $e_1$ is a submersion
(or at least that it is transverse to $\{0\}$).
Let us then denote by $\Sigma_1$ the manifold $e_1^{-1}(0)$,
 and consider the restriction
$\pi_{|\Sigma_1}$ to $\Sigma_1$ of the projection on the second factor.
 This map is $C^1$, and it is Fredholm
of index $0$. Moreover a map $f_0\in F$ is Morse (on $B^n$ ) if and only if 
it is a regular value of $\pi_{|\Sigma_1}$, which means that
the differential of this map is onto at each point $(x,f_0)$
of $\Sigma_1$. These claims  are proved in \cite{AR}, Section 
19, the argument is recalled in Section \ref{ssfredholm} for the convenience of
the reader, see Proposition \ref{lineartrans}. We have proved that 
the set $N\subset F$ of non-Morse functions can be written
$$
N=CV(\pi_{|\Sigma_1}),
$$
where $CV$ denotes the set of critical values. By the theorem 
of  Sard and Smale (see Section \ref{base}),
this set is Baire-meager, and it  has zero measure in $F$, in a sense that will be made precise
in Section \ref{base}.

Let us now present, for $r\geq 3$, a slightly different approach which
has the advantage of avoiding the use of the Theorem of Sard and Smale.
Denoting by $S^n$ the set of symmetric
$d\times d$ matrices,  we start with  the evaluation map
\begin{align*}
e_2:B^n\times F& \lto \Rm^n\times S^n\\
(x,f)&\lmto (df(x),d^2f(x))
\end{align*}
which is $C^{r-2}$ and, when $r \geq 3$ and $F=C^r(\bar B^n,\Rm)$,
is a submersion.
Let us denote  by $\tilde A\subset \Rm^n\times S^n$ the
subset of points $(0,H)$, with $H$ singular.
Note that $\tilde A$ is an algebraic submanifold of codimension 
$n+1$ in $\Rm^n\times S^n$, hence a finite union of smooth submanifolds of codimension
at least $n+1$. The set $N\subset F$ of non-Morse functions
can be written
$$
N=\pi(\Sigma_2), \quad \Sigma_2=e_2^{-1}(\tilde A).
$$
It is best here to first consider that $F$ is finite dimensional (but that $e_2$ is still
a submersion).
Then,  $\Sigma_2$ is a finite union of manifolds of dimension 
less than $\dim F$. This implies that $N=\pi(\Sigma_2)$ is rectifiable of dimension 
less than 
$\dim F$, or in other words it is rectifiable of positive codimension in $F$.
This implies that $N$ has zero measure, but is a much more precise information, which was 
obtained without the use of Sard Theorem.
This reasoning can be extended to the case where $F$ is not finite dimensional
with the help of an appropriate notion of rectifiable sets recalled in Section \ref{base}.
More precisely, we know that $\Sigma_2$ is a finite union of manifolds
of codimension at least $n+1$. 
Since $\pi$ is obviously Fredholm of index 
 $n$,
 we conclude 
by Proposition \ref{direct} that  $\pi(\Sigma_2)$
is rectifiable of codimension $1$.
 We obtain:

\begin{thm}
Let $N\subset C^r(\bar B^n,\Rm)$ be the set of  functions which are not Morse on $B^n$.

If $r\geq 2$, then $N$  is a countable union of closed sets
with empty interior, it has zero measure (in the sense of Haar or Aronszajn).

If $r\geq 3$,
it is  rectifiable of positive codimension. 
\end{thm}

The concepts of sets of zero measure (Haar null or Aronszajn null sets) in 
separable Banach spaces  used in this
statement are recalled in Section \ref{base}, together with the concept of 
rectifiable set of positive codimension. Each point in the statement is the result
of one of the strategies of proof exposed above, notice that none of these statements contains
the other.

As a second illustration, we consider a smooth manifold $A\subset \Rm^m$
and describe the set $NA\subset  F=C^r(\bar B^n,\Rm^m)$ of maps 
which are not transverse to $A$ on $B^n$. 
We consider the evaluation map 
\begin{align*}
E_0:B^n\times F& \lto \Rm^m\\
(x,f)&\lmto f(x).
\end{align*}
This map is $C^{r}$ and, for $r\geq 1$, it is a submersion.
We then have 
$$
NA=CV(\pi_{|\Sigma_0}), \quad \Sigma_0=E_0^{-1}(A),
$$
and $\pi_{|\Sigma_0}$ is $C^{r}$ and Fredholm of index $i=n-c$, where $c$ is the codimension 
of $A$, as follows from Proposition \ref{lineartrans}. If $r\geq n-c+1$, we can apply the theorem 
of  Sard and Smale (see Section \ref{base}), and obtain that 
this set is Baire-meager, and  has zero measure in $F$, in a sense that will be made precise
in Section \ref{base}. When $c>n$, we also conclude that $NA=\pi(\Sigma_0)$ is rectifiable of 
positive codimension.

The second approach, which is useful for $c\leq n$,  consists in using the evaluation map
\begin{align*}
E_1:B^n\times F& \lto \Rm^m\times L(\Rm^n,\Rm^m)\\
(x,f)&\lmto (f(x),df_x).
\end{align*}
This map is $C^{r-1}$, and, for $r\geq 2$, it is a submersion.
Let us denote by $\tilde A$ the set of pairs $(y,l)\in \Rm^m \times L(\Rm^n,\Rm^m)$
such that $y\in A$ and $l(\Rm^n)+T_yA\subsetneq \Rm^m$.
We then have 
$$
NA=\pi(\tilde \Sigma_1), \quad \tilde \Sigma_1=E_1^{-1}(\tilde A).
$$
 We conclude as above that $NA$ is rectifiable of positive codimension
in view of the following Lemma:

\begin{lem}\label{babyEML}
 The set $\tilde A$ is a countable union of smooth manifolds of codimension more than $n$.
\end{lem}

\proof
Locally, there exists a submersion $F:\Rm^m\lto \Rm^c$  such that $A=F^{-1}(0)$.
Then, the set $\tilde A$ is the preimage by the local submersion
$$
\Rm^m\times L(\Rm^n,\Rm^m)  \ni (y,l)\lmto 
(F(y),dF_y\circ l)\in \Rm^c \times L(\Rm^n,\Rm^c)
$$
of the set 
$B:=\{0\}\times L_S(\Rm^n,\Rm^c)
$
where $L_S$ is the set of singular linear maps from $\Rm^n$ to $\Rm^c$
(maps of rank less than $c$).
It is well-known that $L_S(\Rm^n,\Rm^c)$ is an analytic submanifold
of codimension $n-c+1$ in $L(\Rm^n,\Rm^c)$, hence $B$ is an analytic submanifold
of codimension $n+1$. As a consequence, $B$ is a finite union of 
smooth submanifolds of codimension at least $n+1$, hence so is $\tilde A$.
\qed

As a conclusion, we obtain:

\begin{thm}
 Let $A$ be a smooth submanifold of $\Rm^m$ of codimension $c$. 
\begin{itemize}
 \item For $r\geq n-c+1$, The set $NA$ is Baire meager and Aronzajn null (hence Haar null)
in $C^r(\bar B^n,\Rm^m)$.
\item For $r\geq 2$, the set $NA$ is rectifiable of positive codimension in  $C^r(\bar B^n,\Rm^m)$,
it is thus Baire meager and Aronszajn null.
\end{itemize}

\end{thm}

It is worth observing that the second statement contains the first one, except for the case where
$c=n$ and $r=1$.
Our goal in the present paper is to develop an analog of the second strategy 
presented on the examples
above to prove the Thom transversality Theorem in the space of jets.
This idea was suggested by Gromov, in \cite{G}, page 33, and used  in \cite{EM}, 
Section 2.3., where it is reduced to an  appropriate generalization
of Lemma \ref{babyEML} above. 
This Lemma, which is  stated there without proof, is our Conjecture \ref{EML}.
The main  novelty in the present paper consists in giving a full proof of this conjecture
in the analytic case.
We also explain that this strategy, as in the examples
above, leads to a more precise statement of the Thom transversality Theorem
than
 the usual proof based on the Theorem of Sard:

\begin{thm}\label{TT}
 Let $A$ be a smooth submanifold of $J^p( \bar B^n,Y)$ of codimension $c$, where $Y$ is a finite dimensional
 separable  manifold.
For $r\ge p+1$, let $NA\subset C^r(\bar B^n,Y)$ be the set of maps whose $p$-jet is not transverse to $A$.
\begin{itemize}
\item If $c\geq n+1$ and $r\geq p+ 1$, then the set  $NA$  is rectifiable of codimension $c-n$ 
in  $C^r(\bar B^n,Y)$,
it is thus Baire meager and Aronszajn null.
\item If $c\leq n$ and $r\geq p+1+n-c$, the set $NA$ is Baire meager and Aronzajn null (hence Haar null)
in $C^r(\bar B^n,Y)$.
\item 
If $c\leq n$ and $r\geq p+2$ and $A$ is analytic, the set $NA$ is rectifiable of positive codimension 
in  $C^r(\bar B^n,Y)$,
it is thus Baire meager and Aronszajn null.
\end{itemize}
\end{thm}

\proof 
For  completeness, we first quickly recall the usual proof of the 
Thom transversality Theorem, as given in \cite{AR}, which yields the second  point of the Theorem.
We consider the evaluation map (with $F=C^r(\bar B^n,Y)$):
\begin{align*}
E_p:B^n\times F& \lto J^p(\bar B^n, Y)\\
(x,f)&\lmto j^p_xf.
\end{align*}
This map is $C^{r-p}$, and it is a submersion, see \cite{AR},
Theorem 10.4. It follows from Proposition \ref{lineartrans} below that 
$$
NA=CV(\pi_{|\Sigma}), \quad \Sigma=E_p^{-1}(A).
$$
Moreover, the map $\pi_{|\Sigma}$ is Fredholm of index $i=n-c$.
We conclude from the Theorem of Sard and Smale (Theorem \ref{SST} below)
that $NA$ has zero measure and is Baire meager. If, in addition, the codimension of $A$ is 
larger than $n$, then so is the codimension of $\Sigma$, and 
we can conclude directly by the ``Easy Part'' of the theorem of Sard and Smale
that $NA$ is rectifiable of positive codimension.
For the case where $c\leq n$, we obtain the proof of the last point of the theorem 
by considering the evaluation map
\begin{align*}
E_{p+1}:B^n\times F& \lto J^{p+1}(\bar B^n, Y)\\
(x,f)&\lmto j^{p+1}_xf,
\end{align*}
which is a $C^{r-p-1}$ submersion, and the set 
\begin{equation*}
\tilde A\stackrel{\rm def}{=}\big\{ j^{p+1}_xf\in J^{p+1}(\bar B^n, Y): 
j^pf \text{ is not transverse to }A\text{ at }x \big\}\subseteq J^{p+1}(\bar B^n, Y).
\end{equation*}
By definition, we have
$$
NA=\pi(\tilde \Sigma), \quad \tilde \Sigma=E_{p+1}^{-1}(\tilde A).
$$
The last point of Theorem \ref{TT}, without the additional restriction on $A$,
 would be a consequence of the following Conjecture.
The cases of the conjecture that we will be able to prove imply Theorem \ref{TT}.
\qed

\begin{conj}\label{EML}
 If $A$ is a smooth submanifold of $J^p(\bar B^n,Y)$ of codimension $c\leq n$,
then $\tilde A$ is a countable union of smooth submanifolds of codimension 
more than $n$ in $J^{p+1}(\bar B^n,Y)$.
\end{conj}
In view of Proposition \ref{recteasy}, it would even be enough for
our applications to prove that $ \tilde A$ is rectifiable of codimension $n+1$
(in the sense of Section \ref{sec:rect}).
We come back to this conjecture in Section \ref{section:conjecture}, where we obtain some 
special cases, see Theorem \ref{NDT} and \ref{AT}, which  are sufficient 
to imply the third point of Theorem \ref{TT}.
In Section \ref{base}, we recall several mathematical notions which have been 
used in this introduction.

\section{Small sets, rectifiable sets, the theorem of Sard and Smale}\label{base}

\subsection{Some notions of small sets}

Let $F$ be a separable Banach space. We define below three
translation invariant  $\sigma$-ideals of  subsets
of $F$. A  $\sigma$-ideal is a family $\mF$ of subsets of $F$ such that 
\begin{align*}
A\in \mF,  A'\subset A  &\Rightarrow A'\in \mF,\\
\forall n\in \Nm, A_n \in \mF     &\Rightarrow \cup_{n\in \Nm} A_n \in \mF.
\end {align*}

A subset $A\subset F$ is called \textbf{Baire-meager}
if it is contained in a countable union of closed sets with empty interior.
The Baire Theorem  states that a Baire-meager subset
of a Banach space has empty interior.

A subset $A\subset F$ is called \textbf{Haar-null}
if there exists a  Borel  probability measure 
$\mu$ on $F$ such that $\mu(A+f)=0$ for all $f\in F$.
The equality $\mu(A+f)=0$ means that the set $A+f$
is contained in a Borel set $\tilde A_f$ such that 
$\mu(\tilde A_f)=0$.
A countable union of Haar-null sets is Haar-null, 
see \cite{Cr:72,BL} and \cite{prev} for the non-separable case.

A subset $A\subset F$
is called \textbf{Aronszajn-null}
if, for each sequence $f_n$ generating a dense subset of $F$,
there exists a sequence $A_n$ of Borel subsets of $F$
such that $A\subset \cup_n A_n$ and such that,
for each $f\in F$ and for each $n$, the set
$$
\{x\in \Rm: f+xf_n\in A_n\}\subset \Rm
$$
has zero Lebesque measure. 
A countable union of Aronszajn-null sets is Aronszajn-null,
and each Aronszajn-null set is Haar null, see \cite{Ar:76,BL}.

The notion of \textbf{probe} allows a simple criterion for proving that a Borel set $A$
is Haar or Aronszajn null. A probe for $A$ is a finite dimensional vector space  
$E\subset F$ such that $(A+f)\cap E$ has Lebesgue measure zero in $E$ for each $f\in F$.
It is easy to see that $A$ is Haar null if there exists a probe for $A$.
In the sense of Aronszajn, we have (see \cite{Z}, Proposition 4.3):

\begin{lem}\label{probe}
Let $A\subset  F$ be a Borel set.
If the set of probes for $A$ contains a non-empty  open set in the space of finite
dimensional subspaces of $F$, then $A$ is Aronszajn null.
\end{lem}

\proof
Let $f_n$ be a sequence of points of $F$ generating a dense subspace.
Under the hypothesis of the lemma, there exists $N\in \Nm$ and
a probe $E$ such that $E\subset \text{Vect}(f_n,n\leq N)$.
Then, the space $F_N:=\text{Vect}(f_n,n\leq N)$ is itself a probe for $A$.
By standard arguments, (see \cite{BL}, Proposition 6.29 or \cite{Z}), we conclude that 
$A=\cup_{n\leq N} A_n$, where each  $A_n$ is a  Borel set such that  the set
$$
\{x\in \Rm: f+xf_n\in A_n\}\subset \Rm
$$ 
has zero measure for each $f\in F$. Since this holds for each sequence
$f_n$ with dense range, we conclude that $A$ is Aronszajn null.
\qed

If $X$ is a separable manifold modeled on a separable Banach space $F_X$, we 
also have notions of Baire meager, Haar null and Aronszajn null subsets of $X$.
We say that $A\subset X$ is 
Baire meager, Haar null or Aronszajn null if, for each $C^1$ chart
$\varphi: B_X\lto X$, the set $\varphi^{-1}(A)$ is Baire meager, Haar null or Aronszajn null,
where $B_X$ is the open unit ball in $F_X$.
Baire meager sets can also be defined directly as subsets of countable unions of 
closed sets with empty interior in the Baire topological space $X$.

The situation is slightly more problematic with Haar null or Aronszajn-null sets,
because these $\sigma$-ideals are not invariant by $C^1$ diffeomorphisms. 
As a consequence, being Haar null or Aronszajn null in the Banach space $F_X$
seen as a Banach Manifold is a stronger property than being Haar null or Aronszajn null in
$F_X$ seen as a Banach space. This ambiguity in terminology should not cause problems
 in the sequel. Many other notions of sets of zero measure in nonlinear spaces
have been introduced, see for example the survey \cite{HK}.

\subsection{Fredholm maps}\label{ssfredholm}

Given Banach spaces $F$ and $B$,
a continuous linear map $L:F\lto B$ is called Fredholm
if its kernel is finite dimensional and if its range is closed
and has finite codimension.
We say that $L$ is a Fredholm linear map of type $(k,l)$
if $k$ is the dimension of the kernel of $L$
and $l$ is the codimension of its range. The index of $L$
is the integer $k-l$. Recall that the set of Fredholm linear 
maps is open in the space of continuous linear maps (for the norm
topology), and that the index is locally constant, although
the integers $k$ and $l$  are not.
They are lower semi-continuous.
When $F$ and $B$ have finite dimension $n$ and $m$,
then the index of all linear maps is 
$i=n-m$.

The following essentially comes from Section 19 of \cite{AR}.

\begin{prop}\label{lineartrans}
Let $F,X$ be  Banach spaces such that $X$ has finite dimension $n$. Let 
$l:F\times X\lto \Rm^c$
be a surjective continuous  linear map, let $K$ be the kernel of $l$, 
and let $k$ be the restriction to $K$ of the projection 
$(f,x)\lmto f$.

Then $k$ is Fredholm of index $n-c$. Moreover, 
it is onto if and only if
the restriction $l_0$ of $l$ to $\{0\}\times X$ is onto.
\end{prop}

\proof
Let us denote by $X_0$ the space $\{0\}\times X$,  by $F_0$
the space $F\times \{0\}$, and by $K_0$ the intersection
$K\cap X_0$. 
To prove that the continuous linear map $k$ is Fredholm,
we write 
$$
F\times X=K_1\oplus K_0\oplus X_1\oplus F_1
$$
where
\begin{itemize}

\item $F_1\subset F_0$  and $F_1\oplus (K+X_0)=F\times X$.
 Such a space exists because $K+X_0$
has finite codimension, and because $F_0+K+X_0=F\times X$.
\item $K_1\oplus K_0= K$.
\item $X_1\oplus K_0=X_0$.
\end{itemize}
Denoting by $\pi$ the projection on the first factor, we see that
the restriction of $\pi$ to $K_1\oplus F_1$
is an isomorphism onto $F$. 
This implies that the map $k$ is conjugated to linear map:
\begin{align*}
K_1\oplus K_0&\lto K_1\oplus F_1 \\
\kappa_1+\kappa_0&\lmto \kappa_1+0,
\end{align*}
which is Fredholm of index $i=\dim K_0-\dim F_1$.
We obtain that 
$$i=(\dim K_0+\dim X_1)-(\dim X_1+\dim F_1)=n-c. 
$$
The linear Fredholm map 
$k$ is onto if and only if its kernel $K_0$ has dimension 
$k=n-c$. On the other hand, the space $K_0$ is also the kernel of $l_0$,
hence it   
has dimension $n-c$ if and only if the  $l_0$ is onto
( $X_0$ has dimension $n$).
We have proved the second part of the statement.
\qed

We now recall a  standard Lemma of differential calculus.
\begin{lem}\label{normalform}
 Let $f:X\lto Y$ be a $C^1$ map and $x_0$ be a point such that 
$df_{x_0}$ has a closed and split 
range $I\subset F_Y$ and a split kernel $K\subset F_X$.
Let $G$ be a supplement of $I$ in $F_Y$.
Then, for each local diffeomorphism 
$\phi:(Y,f(x_0))\lto (I\times G,0)$
there exists a local diffeomorphism 
$\varphi: (I\times K,0)\lto (X, x_0)$ such that 
$$
\phi \circ f\circ \varphi (x_i,x_k)=(x_i,\psi(x_i,x_k))
$$
for some $C^1$ local map
$\psi:I\times K\lto G$.
\end{lem}

This Lemma can be applied in particular to Fredholm maps.

\proof
Let $E$ be a supplement of $K$ in $F_X$.
By considering first an arbitrary local chart 
$\tilde \varphi:(E\times K,0)\lto (X,x_0)$, we write
$$
\phi\circ f \circ \tilde \varphi: (x_e,x_k)\lmto (f_i(x_e,x_k),f_g(x_e,x_k)).
$$
It follows from the definition of $G$ and $I$ that $\partial _{x_e}f_i$
is an isomorphism, hence the mapping
$$
(x_e,x_k)\lmto (f_i(x_e,x_k),x_k)
$$
is a local diffeomorphism   between $(B_X,0)$ and $(I\times K,0)$.
Denoting by $\hat \varphi(x_i,x_k)=(\hat \varphi_e(x_i,x_k),x_k)$
its inverse, we see that 
$$
\phi\circ f \circ \tilde \varphi\circ \hat \varphi(x_i,x_k)=(x_i,\psi(x_i,x_k))
$$
with $\psi(x_i,x_k)=f_g(\hat \varphi_e(x_i,x_k),x_k)$.
\qed

We also recall the constant rank (or rather constant corank) theorem.

\begin{lem}
 Let $f:X\lto Y$ be a $C^1$ map. Assume that there exists an integer $c$ such
 that, for each $x\in X$, 
$df_{x}$ has a closed 
range $I\subset F_Y$ of codimension $c$ and a split kernel $K\subset F_X$.
Let $G$ be a supplement of $I$ in $F_Y$.
Then, near each point $x_0\in X$ there exists a local   diffeomorphism 
$\phi:(Y,f(x_0))\lto (I\times G,0)$
and a local diffeomorphism 
$\varphi: (I\times K,0)\lto (X, x_0)$ such that 
$$
\phi \circ f\circ \varphi (x_i,x_k)=(x_i,0).
$$
\end{lem}

\proof 
We first apply Lemma \ref{normalform} and find charts $\tilde \phi$
and $\varphi$ such that 
$
\tilde \phi \circ f\circ \varphi (x_i,x_k)=(x_i,\psi(x_i,x_k)).
$
The differential of this map has corank $c$ (which is the dimension of $G$)
if and only if $ \partial_{x_k}\psi=0$. We conclude that $\psi$ does 
not depend on $x_k$.
We now set $\hat \phi(x_i,x_g)=(x_i,x_g-\psi(x_i))$, and observe that 
$
\hat \phi \circ\tilde \phi \circ f\circ \varphi (x_i,x_k)=(x_i,0).
$
\qed

\subsection{Rectifiable sets in Banach manifolds}\label{sec:rect}

We use here the definition of rectifiable sets of finite codimension
given in \cite{ARMA}, which extrapolates on \cite{Z}. 
Our terminology, however, differs from that of \cite{ARMA} : 
we call rectifiable here
what we called countably rectifiable there.

The subset $A$ in the Banach space $F$ is a 
\textbf{Lipschitz graph of codimension d} 
if there exists a splitting $F=E\oplus G$, with $\dim G=d$ and a Lipschitz
map $g:E\lto G$ such that 
$$
A\subset \{ x+g(x), x\in E\}.
$$

Let $X$ be a separable manifold modeled on the separable Banach space $F_X$.
A subset $A\subset X$ is \textbf{rectifiable of codimension d} if it is 
a countable union $A=\cup _n \varphi_n(A_n)$ where 
\begin{itemize}
 \item $\varphi_n:U_n\lto X$ is a Fredholm map
\footnote{A Fredholm map of index $i$ between separable Banach manifolds is a $C^1$ map such that the differential is Fredholm of index $i$ at every point (recall that the index is locally constant).}
  of index $i_n$ defined
on an open subset $U_n$ in a separable Banach space $F_n$.
\item $A_n\subset U_n$ is a Lipschitz graph of codimension $d+i_n$ in $F_n$.
\end{itemize}

Note that, by definition, if $A\subset X$ is rectifiable of codimension $d$ then it is rectifiable of codimension $d'$ for all $0\le d'\le d$.
The following properties are proved in \cite{Z} or \cite{ARMA}.

\begin{prop}\label{smallness}
A rectifiable set of positive codimension is Baire meager.
More precisely, it is contained in a countable union 
of closed sets of  positive codimension. 
It is also Aronszanjn null, hence Haar null.
\end{prop}

\begin{prop}\label{direct}
Let $X$ and $Y$ be separable Banach manifolds, and let 
 $f:X\lto Y$ be $C^1$ Fredholm of index $i$, and let $A\subset X$ be rectifiable of codimension 
$d\geq i+1$,
then the direct image  $f(A)$ is rectifiable of codimension $d-i$.
\end{prop}

The following property is almost taken from \cite{ARMA}.

\begin{prop}\label{preimage}
Let $X$ and $Y$ be separable Banach manifolds, let $A\subset Y$
be rectifiable of codimension $d$ and let 
 $f:X\lto Y$ be a $C^1$ map such that, at each point of $f^{-1}(A)$, the differential
$df$ has the following properties with some integer $k\leq d-1$:
\begin{itemize}
 \item  It has a split kernel.
\item It has a closed image of codimension at most $k$.
\end{itemize}
Then, $f^{-1}(A)$ is rectifiable of codimension $d-k$.
In particular, if $f$ is a submersion, then $f^{-1}(A)$ is rectifiable of codimension $d$.
\end{prop}
\proof
In view of Lemma \ref{normalform}, it is enough to prove the statement for 
maps of the form $(x_i,x_k)\lmto (x_i,\psi(x_i,x_k))$, 
where $(x_i,x_k)\in I\times K$, and $\psi$ takes value in $G$, a supplement
of $I$ in $B_Y$ (hence $\dim G\leq k$).
We also consider that $A\subset  I\times G$ is rectifiable of codimension $d$.
Then the projection $A_I$ of $A$ on $I$ is rectifiable of codimension 
$d-k$ in $I$. In view of the special form of the map we consider,
the preimage of $A$ is contained in $A_I\times K$, which is rectifiable 
of codimension $d-k$ in $I\times K$ because $A_I$ is rectifiable of codimension
$d-k$ in $I$.
\qed

We express the following results in the context of Banach spaces to avoid some technical complications. If $F$ is a separable Banach space, then we define the separable Banach spaces
$C^p(\bar B^n,F)$ as in the introduction.
The following result was proved in \cite{ARMA}:

\begin{prop}
 If $A\subset C^p(\bar B^n,F)$ is rectifiable of codimension $d$,
 and $p'\geq p$, then $A\cap C^{p'}(\bar B^n,F)$
 is rectifiable of codimension $d$ in $C^{p'}(\bar B^n,F)$.
\end{prop}

This results allows to define sets of positive codimension 
in the Frechet space $C^{\infty}(\bar B^n,F)$, see \cite{AR}.
The following result makes precise the simple fact that
``most'' $n$-parameter families avoid sets of codimension $d$
when $d>n$. 

\begin{prop}
 Let $F$ be a separable Banach space, and $A\subset F$
 a rectifiable set of codimension $d$.
 For $n<d$, The set $\mA\subset C^1(\bar B^n,F)$ of maps $f$ such that 
 $f(B^n)\cap A \neq \emptyset$ is rectifiable of codimension 
 $d-n$.
\end{prop}

\proof
This is just a variant of the methods of proof used in the introduction.
We consider the evaluation map
$$
E_0: B^n\times C^p(\bar B^n,F)\lto F,
$$
which is a $C^1$ submersion.
We conclude from Proposition \ref{preimage}
that $E_0^{-1}(A)$ is rectifiable of codimension $d$ in
$B^n\times C^p(\bar B^n,F)$.
The set $\mA$, which is the projection of $E_0^{-1}(A)$
on the second factor, is thus rectifiable of codimension $d-n$.
\qed
The ``easy case'' of the transversality theorem also has a natural  analog
in terms of rectifiable sets. 

\begin{prop}\label{recteasy}
 Let $F$ be a separable Banach space, and $A\subset J^p(\bar B^n,F)$
 a rectifiable set of codimension $d$.
 For $n<d$, The set $\mA\subset C^{p+1}(\bar B^n,F)$ of maps $f$ such that 
 $j^pf(B^n)\cap A \neq \emptyset$ is rectifiable of codimension 
 $d-n$.
\end{prop}

\proof
It is the same as above, using the evaluation map 
$E_p:(x,f)\lmto j^p_xf$.
\qed
\subsection{The Theorem of Sard and Smale}

\begin{thm}
Let $X$ be a smooth (separable) manifold of dimension $n$,
 $Y$ a smooth manifold of dimension $m$, and 
let $f:X\lto Y$ be a $C^r$ map. 
If $r\geq 1+(n-m)$, and  $m\leq n$, then the set $CV(f)$ of
critical values of $f$ has zero measure in $Y$.
\end{thm}

The theorem also holds in the case where $r\geq 1$ and $n\leq m-1$,
which is sometimes called the easy case of Sard's theorem.
In this case, however, the set $CV(f)$ is just the image $f(X)$, which
is rectifiable of dimension $n$ in $Y$. Since $n<m$, this implies the result, but is 
a much finer information. The theorem of Sard was extended by Smale to the
 infinite dimensional case. We give below a more precise statement:

\begin{thm}\label{SST}
Let $X$ and $Y$  be  separable smooth  manifolds modelled on separable Banach spaces,
and 
let $f:X\lto Y$ be a $C^r$ Fredholm map of index $i$. 
\begin{itemize}
\item If $i\geq 0$ and $r\geq 1+i$,  then the set $CV(f)$ of
critical values of $f$ is Aronszajn null (hence Haar null) and Baire meager  in $Y$.
\item If $i<0$ and $r\geq 1$, then the set $f(X)=CV(f)$ is rectifiable of codimension $-i$
in $Y$  (It is thus Aronszajn-null  and Baire meager).
\end{itemize}
\end{thm}

\proof
The second part of the statement (the ``easy case''),
 is a special case  of Proposition \ref{direct}.
 Let us focus on the first part.
 Let $P$ be the set of critical points of $f$, so that $CV(f)=f(P)$.
 
 We claim that each point $x_0$ of $P$ has a closed neighborhood
 $\tilde P$ in $P$  such that $f(\tilde P)$ is closed and Aronszajn null. 
 Since $P$ is a separable metric space, it has the Lindel\"of property,
 and  it  can be covered by countably
 many such local sets  $\tilde P$.
 As a consequence, the claim implies the statement.
 
Since the claim is local, we identify $X$ with its Banach model $B_X$ and
$x_0$ with $0$, and similarly $(Y,f(x_0))$ with $(B_Y,0)$.
 Let $I$ be the range of $df_0$, and let 
 $G$ be a supplement of 
 $I$ in $F_Y$, note that $G$ has finite dimension $l$.
 By Lemma \ref{normalform}, there exists a local diffeomorphism
 $\varphi: (I\times K,0)\lto (B_X,0)$ such that 
 $$
 f\circ \varphi (x_i,x_k)=x_i+\psi(x_i,x_k),
 $$
 where $\psi:I\times K\lto G$ is $C^{r}$.
 We conclude that 
 $$
 \phi(f(\tilde P))\subset CV(\phi\circ f\circ \varphi)
 \subset \bigcup_{x_i\in I}\big( x_i+CV(\psi_{x_i})\big)
 $$
  where $\psi_{x_i}:K\lto G$ is the map $x_k\lmto \psi(x_i,x_k)$.
 In view of the finite dimensional Sard theorem
 (applied to  $\psi_{x_i}$), we conclude  that $G$ 
 is a probe for $\phi(f(\tilde P))$. Since it admits a probe,
 this set is Haar-null. Moreover, since the set of supplements
 of the $I$ is open in the space of $l$-dimensional subspaces of 
 $F_Y$, we conclude from Lemma \ref{probe} that $\phi(f(\tilde P))$
 is Aronszajn null.
 Since this holds for each local chart $\phi$, we have proved that
 $f(\tilde P)$ is Aronszajn null in $Y$.
 
 Finally, let us prove that $f(\tilde P)$ is closed, 
 or equivalently that $\phi(f(\tilde P))$ is closed, provided $\tilde P$
 is chosen bounded and closed in $B_X$.
 Let $x^n$ be a sequence in $\tilde P$, such that $f(x^n)$ has a limit 
 $y^{\infty}$, we have to prove that $y^{\infty}\subset f(\tilde P)$.
 Let us denote by $(x_i^n,x_k^n)$ the sequence $\varphi^{-1}(x^n)$.
 Since $K$ is finite dimensional and $x^n$ is bounded,
 we can assume by taking a subsequence that $x_k^n$ has a limit $x_k^{\infty}$.
 On the other hand, since $f(x^n)\lto y^{\infty}$,
 we conclude that 
 $$
 (x_i^n,\psi(x_i^n,x_k^n))=\phi\circ f\circ \varphi(x^n_i,x^n_k)=
 \phi\circ f(x^n)\lto \phi(y^{\infty}),
 $$
 hence $x_i^n$ has a limit $ x_i^{\infty}$ (which is the first component of 
 $\phi(y^{\infty}))$.
 The sequence  $(x_i^n,x_k^n)$ is thus convergent, hence so is 
 $x_n=\varphi(x_i^n,x_k^n)$. Since $\tilde P$ is closed, the limit
 $x^{\infty}$ belongs to $\tilde P$, and $y^{\infty}=f(x^{\infty})$.
 \qed

\section{Some cases of the conjecture}\label{section:conjecture}

In this section, we consider two finite dimensional smooth manifolds 
$X$ and $Y$. For $p\in \Nm$, we denote by $J^p(X,Y)$ the space of $p$-jets
of functions $X\lto Y$. 
See for example \cite{H,GG} for some details on jet bundles.
Note that $J^0(X,Y)=X\times Y$, and it will also be convenient to
consider that $J^{-1}(X,Y)=X$.
For $p\leq p'$ we have a natural projection
$$
\pi^{p'}_{p} :J^{p'}(X,Y)\lto J^p(X,Y).
$$
When $p=-1$, this is just the source map $j^{p'}_x f\lmto x$.
When $p\in \Nm$, the bundle 
$$
\pi^{p+1}_{p} :J^{p+1}(X,Y)\lto J^p(X,Y)
$$
has a natural affine structure, we denote by $\mF_p^{p+1}(a)$ the fiber $(\pi_p^{p+1})^{-1}(a)$,
for $a\in J^p(X,Y)$.
Given a submanifold  $A\subseteq J^p(X,Y)$ 
 of class $C^r, r\ge 1$, we define 
\begin{equation}
\label{eq:tilde}
\tilde A\stackrel{\rm def}{=}\big\{ j^{p+1}_xf\in J^{p+1}(X,Y): 
j^pf \text{ is not transverse to }A\text{ at }x \big\}\subseteq J^{p+1}(X,Y).
\end{equation}
If $j^p_x f$ is an element of $J^p(X,Y)$ and $0\le k\le p-1$, we define the vector subspace $E^{k}(j^p_xf)$ by
\begin{equation}
\label{eq:definition E}
E^{k}(j^p_xf)\stackrel{\rm def}{=}d(j^{k}_x f)(T_x X)\subseteq T_{j^{k}_x f}(J^{k}(X,Y))
\end{equation}
Here $d(j^{k}_x f)$ is the tangent map at $x$ of $j^{k}f\colon X\to J^{k}(X,Y)$.
 Note that the subspace $E^k(j^p_xf)$ depends just on $j^{k+1}_xf$ 
and that its dimension is always equal to $\dim X$. 
We also extend the definition to $k=-1$ in a trivial way by setting $E^{-1}(j^p_xf)=T_x X$. 
We have 
$$
\tilde A=
\big\{ z\in J^{p+1}(X,Y): 
T_{\pi^{p+1}_p z} A+E^p(z)
\subsetneq T_{\pi^{p+1}_p z }J^p(X,Y)
\big\}.
$$

\begin{conj}\label{EML}
 The set $\tilde A$ is a countable union of submanifolds
of codimension more than $n=\dim X$.
\end{conj}

This conjecture is stated as a Lemma in \cite{EM}, 
but not proved.
The statement of the theorem is obvious when the codimension $c$ of $A$ is larger
than $n$, we assume from now on that $c\leq n$.
 We will use the notation
\[
\tilde A_a=\tilde A\cap \mF_p^{p+1}(a). 
\]
We say that the point  $a\in A$ is degenerate if 
\begin{equation}
\label{eq:deg}
\pi^p_{p-1}(T_a A)+ E^{p-1}(a)\subsetneq   T_{\pi^p_{p-1}a}J^{p-1}(X,Y).
\end{equation}
(for $p=0$ the map $\pi^0_{-1}$ is the projection from $J^0(X,Y)$ to $X$). 
If $a$ is degenerate, then $\tilde A_a=\mF^{p+1}_p(a)$.
The manifold $A$ can be decomposed as the disjoint unions
$A=A_0\cup A_1$, where $A_0$ is the set of degenerate points of $A$ and $A_1$ is the 
set of non-degenerate points (by definition, the other points). The set $A_1$
is an open submanifold of $A$, hence it is itself a submanifold of $J^p(X,Y)$,
and  
\begin{equation}
\label{eq:disjoint union}
\tilde A=(\pi^{p+1}_p)^{-1}(A_0) \cup \tilde{A_1}.
\end{equation}
Let us first treat the special case where $A=A_1$ (we then say that $A$ is non-degenerate).
 Note that this condition holds for example if $A$ is transverse to 
the fibers of the projection $\pi^p_{p-1}$. This  condition also holds when $p=0$.
The following result also implies Lemma \ref{babyEML}:

\begin{thm}\label{NDT}
 Let $A$ be a non-degenerate $C^r$ submanifold in $J^p(X,Y)$.
Then $\tilde A$ is a countable union of $C^{r-1}$ submanifolds of codimension 
larger than $n=\dim X$ in $J^{p+1}(X,Y)$. 
\end{thm}

This result implies that $\tilde A_1$ is a countable union of submanifolds of codimension at least $n+1$.
 In order to prove the 
conjecture, we would also need to prove that the manifold $A_0$ has codimension 
$n+1$. We are not able to prove this statement in the general case,
hence we will restrict to the analytic case.
We say that the submanifold $A\subseteq J^p(X,Y)$ is \emph{analytic}
 if for every $a=j^p_x f\in A$ there exist charts $\psi_X$ and $\psi_Y$ on $X$ and $Y$, 
respectively defined on a neighborhood of $x$ in $X$ and a neighborhood of $f(x)$ in $Y$, 
such that the induced chart $\psi$ on $J^p(X,Y)$, defined on a neighborhood $U_a$ of $a$,
 makes $A$ analytic, i.e.
\[
\psi(A\cap U_a)=\bigcap_i F_i^{-1}(0)
\]
for a finite family of analytic functions $F_i\colon\psi(U_a)\to\R$.
When $A$ is analytic, we manage to study $A_0$ by recurrence using Theorem \ref{NDT},
and obtain:

\begin{thm}\label{AT}
 Let $A$ be an analytic  submanifold in $J^p(X,Y)$.
Then $\tilde A$ is a countable union of $C^{r-1}$ submanifolds of codimension 
larger than $n=\dim X$ in $J^{p+1}(X,Y)$. 
\end{thm}

\subsection{The non-degenerate case}\label{NDproof}

We assume here that the $C^r$ manifold $A\subset J^p(X,Y)$ is non-degenerate,
which means that 
\begin{equation}
\label{eq:nondeg}
\pi^p_{p-1}(T_a A)+ E^{p-1}(a)=  T_{\pi^p_{p-1}a}J^{p-1}(X,Y)
\end{equation}
for each $a\in A$, and prove Theorem \ref{NDT}.
To study the set $\tilde A_a$, we define, more generally, the set
\[
Z_{a,V}=\left\{ \hat a\in\mF_p^{p+1}(a):
V+E^{p}(\hat a)\subsetneq T_a J^p(X,Y) \right\}
\subseteq \mF_p^{p+1}(a)
\] 
associated to a point $a\in J^p(X,Y)$ and a subspace $V\subset T_aJ^p(X,Y)$.
Then, we have
$$
\tilde A_a=Z_{a,T_aA}.
$$
We decompose $Z_{a,V}$  as
\[
Z_{a,V}=\bigcup_{r= \dim{V}}^{\dim J^p(X,Y)-1} Z_{a,V}^r 
\]
where
\[
Z_{a,V}^r\stackrel{\rm def}{=}
\left\{ \hat a\in(\mF^{p+1}_p)(a):\dim\Big(V+E^{p}(\hat a)\Big)=r \right\}.
\]
This decomposition obviously yields a decomposition
$
\tilde A=\cup \tilde A^r$, where 
$$
\tilde A^r_a:=
 \tilde A^r\cap \mF^{p+1}_p(a)
=Z^r_{a,T_aA}.
$$
The following result implies that $\tilde A^r$ is a $C^{r-1}$  submanifold
of codimension at least $n+1-c$ in $(\pi_p^{p+1})^{-1}(A)$, 
hence a submanifold of codimension at least $n+1$ in $J^{p+1}(X,Y)$,
which proves Theorem \ref{NDT}.

\begin{prop}
\label{prop:core proposition}
Let $a\in J^p(X,Y)$ and $V$ be a vector subspace of $T_a J^p(X,Y)$ of dimension $m$ and codimension $c\ge 1$ such that
\begin{equation}
\label{eq:projected transversality}
\pi^p_{p-1}(V)+E^{p-1}(a)=T_{\pi^p_{p-1}a}J^{p-1}(X,Y).
\end{equation}
We have:
\[
\codim_{\mF^{p+1}_p(a)}\,Z_{a,V}\ge n+1-c.
\]
%More precisely, the set $Z^r_{a,V}$ can be written locally as the preimage $F_{a,V}^{-1}(0)$ of an algebraic submersion 
More precisely, the set $Z^r_{a,V}$ is locally contained in the preimage
 $F_{a,V}^{-1}(0)$ of an algebraic submersion 
$$
F_{a,V}:\mF_p^{p+1}(a)\lto \Rm^ {\theta}
$$
whose coefficients depend smoothly on $(a,V)$,
with $\theta=n+1-c+ 2(\dim J^p(X,Y)-1-r).$
\end{prop}

\proof 
Since the result is of local nature, we can suppose without loss of generality that the jet bundles are trivialized. Hence,
\[
J^{p+1}(X,Y)=J^p(X,Y)\times \F^{p+1}_p
\]
where $\F^{p+1}_p$ is the fiber of the projection $\pi^{p+1}_p\colon J^{p+1}(X,Y)\to J^p(X,Y)$. 
We have thus the identification 
%$(\pi_p^{p+1})^{-1}(a)\cong \F^{p+1}_p$
$\mF_p^{p+1}(a)\cong \F^{p+1}_p$.
 Hence the sets $Z^r_{a,V}$ can be regarded as subsets of $\F^{p+1}_p$: denoting by $z$ the elements of $\F^{p+1}_p$, we have
\[
Z^r_{a,V}=\left\{ z\in\F^{p+1}_p: \dim \Big(V+E^{p}(a,z)\Big)=r  \right\}.
\]
We assume from now on that $m\le r\le \dim J^p(X,Y)-1$. Let us pick, for any $z$ in $\F^{p+1}_p$, a function $f_z$ such that
\begin{equation*}
j^{p+1}_x f_z=(a,z).
\end{equation*}
Here $x\in X$ is the base-point of $a$.
 Let us also choose a basis $v_1,\dots,v_m$ of $V$.
 Let us call $M_{a,V}(z)$ (or just $M(z)$)  the matrix whose columns are, in the order, the following vectors (belonging to $T_a J^p(X,Y)$)
\begin{equation}
\label{eq:columns}
v_1,\dots,v_m,\,\partial_{x_1}j^p_x f_z,\dots,\partial_{x_n}j^p_x f_z
\end{equation}
expressed in a convenient basis of $T_a J^p(X,Y)$ 
 which we shall explicit shortly.
  Note that the vectors $\partial_{x_j}j^p_x f_z, 1\le j\le n$ form a basis of $E^{p}(a,z)$. It is then clear that
\[
z\in Z^r_{a,V}\Leftrightarrow \rank M_{a,V}(z)=r,
\]
or equivalently
\begin{align*}
z\in Z^r_{a,V}\Leftrightarrow\begin{cases} \deter N(z)=0\qquad\forall\ \text{square submatrix $N$ of $M_{a,V}$ of size $r+1$}
\\
\deter N(z)\neq 0 \qquad\text{ for some square submatrix $N$ of size $r$}.
\end{cases}
\end{align*}
We will now study more precisely these equations with the help of 
an appropriate system of  local coordinates.
Locally, we have the identifications
\begin{equation}
\label{eq:trivializations}
\begin{aligned}
J^p(X,Y)&=J^{p-1}(X,Y)\times\F^p_{p-1}
\\
J^{p+1}(X,Y)&=J^{p-1}(X,Y)\times\F^p_{p-1}\times\F^{p+1}_p,
\end{aligned}
\end{equation}
and both $\F^{p+1}_p$ and $\F^p_{p-1}$ can be identified with real vector spaces. 
More precisely, we fix once for all local coordinates $x_1,\dots,x_n$ and $y^1,\dots,y^q$ on $X$
 and $Y$ respectively, this induces the identification $\F^{p+1}_{p}\cong\R^{\dim\F^{p+1}_{p}}=\R^{q\binom{n+p}{n-1}}$ via the isomorphism
\begin{align*}
\big(y^s_\alpha\big)_{1\le s\le q,\,|\alpha|=p+1}\colon \F^{p+1}_{p}&\to \R^{q\binom{n+p}{n-1}}
\\
j^p_x f &\mapsto \big(\partial_{\alpha} f^s(x)\big)_{1\le s\le q,\,|\alpha|=p+1}.
\end{align*}
Here $f^s=y^s\circ f$ is the $s$-th component of $f$, 
$\alpha=(\alpha_1,\dots,\alpha_n)$ is a multi-index in $\N^n$ of length $|\alpha|=p+1$ 
and $\partial_\alpha=\partial^{\alpha_1}_{x_1}\dots\partial^{\alpha_n}_{x_n}$ stands for 
the associated partial derivative. Note that for the isomorphism to be rigorously defined,
 one should specify an order on the set of the involved couples $(s,\alpha)$.
 Since this order will not play any role in the sequel, we do not specify it.
\bigbreak

Concerning $\F^p_{p-1}$, we have the analogous identification $\F^{p}_{p-1}\cong\R^{\dim\F^{p}_{p-1}}=\R^{q\binom{n+p-1}{n-1}}$ via the isomorphism
\begin{align*}
\big(y^s_\alpha\big)_{1\le s\le q,\,|\alpha|=p}\colon \F^p_{p-1}&\to \R^{q\binom{n+p-1}{n-1}}
\\
j^{p}_x f &\mapsto \big(\partial_{\alpha} f^s(x)\big)_{1\le s\le q,\,|\alpha|=p}
\end{align*}
Here the order on the couples $(s,\alpha)$ will play an important role. For reasons which will become clear 
%in Part 4
 in the Lemma \ref{lexicographic property},
 we adopt the following lexicographic order:\footnote{What really matters for our purposes is the lexicographic order with respect to $\alpha$ at fixed $s$. There are several orders satisfying this condition, but for the sake of definiteness we adopt the one described in \eqref{eq:order}.} if $s,s'\in\{1,\dots,q\}$ and $\alpha=(\alpha_1,\dots,\alpha_n)$, $\alpha'=(\alpha'_1,\dots,\alpha'_n)$ are multi-indexes of length $p$, the variable $y^s_\alpha$ strictly precedes the variable $y^{s'}_{\alpha'}$ if and only if
\begin{equation}
\label{eq:order}
s> s'\qquad\text{ or }\qquad \left(s= s'\quad\text{ and }\quad\exists\ k\ge 1:\begin{cases}
\alpha_h=\alpha'_h\quad & \text{for all $1\le h\le k-1$}
\\
\alpha_k >\alpha'_k& 
\end{cases}\right)
\end{equation}
\bigbreak

Summing up the above paragraphs, we will regard 
$\big(y^s_\alpha\big)_{1\le s\le q,\,|\alpha|=p+1}$ and
 $\big(y^s_\alpha\big)_{1\le s\le q,\,|\alpha|=p}$ as coordinates
 respectively on $\F^{p+1}_p$ and $\F^p_{p-1}$, compatible with the
 affine structure of these spaces. The coordinates on $\F^p_{p-1}$ are ordered according 
to \eqref{eq:order}.  Since we denoted by $z$ the elements of $\F^{p+1}_p$, we have
\[
z=\big(y^s_\alpha\big)_{1\le s\le q,\,|\alpha|=p+1}\in\F^{p+1}_p\cong\R^{\dim \F^{p+1}_p}.
\]

\bigbreak
Let us now write more explicitly the vectors $\partial_{x_j}j^p_x f_z$ appearing in \eqref{eq:columns}. According to the decomposition \eqref{eq:trivializations}, an arbitrary $p$-jet $j^p_x f$ writes
\begin{equation}
\label{eq:p-jet}
j^p_x f=\Big(j^{p-1}_x f,(\partial_{\alpha}f^s(x))_{1\le s\le q,\,|\alpha|=p}\Big)\quad\in J^{p-1}(X,Y)\times \F^p_{p-1}.
\end{equation}
The vectors $\partial_{x_j}j^p_x f_z$ are elements of the vector space 
$T_a J^p(X,Y)=T_{\pi^p_{p-1}a}J^{p-1}(X,Y)\times \F^p_{p-1}$.
 By taking the derivative of \eqref{eq:p-jet} with respect to $x_j$ we get:
\begin{align*}
\partial_{x_j} j^p_x f_z&= \Big(\partial_{x_j} j^{p-1}_x f_z,\big(\partial_{\alpha+\delta_j}f^s(x)\big)_{1\le s\le q,\,|\alpha|=p}\Big)
\\
&=\Big(\partial_{x_j} j^{p-1}_x f_z,\big(y^s_{\alpha+\delta_j}\big)_{1\le s\le q,\,|\alpha|=p}\Big)
\end{align*}
where $\delta_j$ is the multi-index
% $\big(0,\dots,0,\underbracket[0.3pt][2pt]{1}_{j\text{-th}},0,\dots,0\big)\in\N^n$. 
%
$\big(0,\dots,0,\underbrace{1}_{j\text{-th}},0,\dots,0\big)\in\N^n$.
Note that the component $\partial_{x_j} j^{p-1}_x f_z$ depends just
 on partial derivatives of order one of $j^{p-1} f_z$ at $x$, i.e.\ it 
depends just on $j^p_x f_z=a$. This justifies the following notation:
\[
e_j(a)\stackrel{\rm def}{=}\partial_{x_j}j^{p-1}_x f_z,\qquad 1\le j\le n.
\]
Note that
\[
\Span \{e_1(a),\dots,e_n(a)\}=E^{p-1}(a).
\]
We have:
\begin{equation}
\label{eq:partial derivative}
\partial_{x_j} j^p_x f_z=\Big(e_j(a),\big(y^s_{\alpha+\delta_j}\big)_{1\le s\le q,\,|\alpha|=p}\Big).
\end{equation}

We can now write the matrix $M_{a,V}(z)$ in the base of 
$T_a J^p(X,Y)=T_{\pi^p_{p-1}a}J^{p-1}(X,Y)\times \F^p_{p-1}$ obtained as the juxtaposition
of an arbitrary base of $T_{\pi^p_{p-1}a}J^{p-1}(X,Y)$ and of the base 
$\big(y^s_\alpha\big)_{1\le s\le q,\,|\alpha|=p}$
 of $\F^p_{p-1}$ in the order which has been specified before:
\[
M_{a,V}(z)= 
\begin{array}{c@{}c}
\left[
  \begin{BMAT}[5pt]{c|c}{c}
    \begin{BMAT}[5pt]{c:c:c}{ccccccc}
      & & \\
      & & \\
      & & \\
      v_1 & \dots & v_m \\
      & & \\
      & & \\
      & &
    \end{BMAT}
    &
    \begin{BMAT}{c}{c|c}
      \begin{BMAT}[5pt]{c:c:c}{ccc}
        & & \\
        e_1(a) & \dots & e_n(a) \\
        & &
      \end{BMAT}
      \\
      \begin{BMAT}[10pt]{c}{c}
        B(z)
      \end{BMAT}
    \end{BMAT}
  \end{BMAT} 
\right] 
& 
\begin{array}{l}
  \\[-10mm] \rdelim\}{4}{40mm}[{\footnotesize $\dim J^{p-1}(X,Y)$}] \\ \\ \\[5mm]  
\rdelim\}{3}{40mm}[{\footnotesize $\dim \mF^p_{p-1}$}] \\ \\
\end{array} \\[-1ex]
%\!\hexbrace{2.6cm}{m=\dim V}\ \ \hexbrace{3.8cm}{n=\dim X}
\end{array}
\]
Here the vectors $v_1,\dots,v_m$ are a basis of $V$ and the vectors $e_1(a),\dots,e_n(a)$,
which  have been defined before, form a base of $E^{p-1}(a)$. Finally, the block $B(z)$ depends just on $z$ and is given by
\begin{equation}
\label{eq:B(z)}
B(z)=
\begin{blockarray}{[ccccc]c}
\vdots & \vdots & \vdots & \vdots & \vdots & 
\\
\ y^{s}_{\alpha+\delta_1} \ &\ y^{s}_{\alpha+\delta_2} \ & \  \cdots\ &\ y^{s}_{\alpha+\delta_{n-1}}\ &\ y^{s}_{\alpha+\delta_n} \ &\ \text{$\leftarrow$ row corresponding to $y^s_\alpha$}
\\ 
\vdots & \vdots & \vdots & \vdots & \vdots & 
\end{blockarray}
\end{equation}
where the rows are ordered according to \eqref{eq:order}. 

For later use, note that the first $m$ columns of $M(z)$ are linearly independent, 
as well as the first $\dim J^{p-1}(X,Y)$ rows, as follows from the hypothesis
 \eqref{eq:projected transversality}.
Indeed, the first $m$ columns are clearly independent because they represent a basis of $V$. The fact that the first $\dim J^{p-1}(X,Y)$-rows are independent is equivalent to the assumption $\pi^p_{p-1}(V)+E^{p-1}(a)=T_{\pi^p_{p-1}a}J^{p-1}(X,Y)$.

\textbf{An intermede of linear algebra.}
We prove, for an arbitrary matrix, the existence of a non-singular
 square submatrix of maximum rank satisfying some special conditions.

Let $M$ be an arbitrary $m\times n$-matrix with real entries. Only in this 
%Part 3
 intermede, 
 $m$ and $n$ are arbitrary integers $\ge 1$, with no relation with the values assumed by the same symbols in the rest of the paper.

Let us establish some notations for submatrices of $M$.  The rows of $M$ are $\Rr(M)=\{1,2,\dots,m\}$ and its columns are $\Cc(M)=\{1,2,\dots,n\}$. A submatrix $N$ of $M$ is determined by its rows 
\[
\Rr(N)\subseteq \{1,2,\dots,m\}
\]
and its columns
\[
\Cc(N)\subseteq \{1,2,\dots,n\}.
\]
We denote $|\Rr(N)|$ and $|\Cc(N)|$ their cardinality. We also denote $i_1(N),i_2(N),\ldots i_{|\Rr(N)|}(N)$ the elements of $\Rr(N)$, and we always assume that the indexes are chosen in such a way that
\[
i_1(N)<i_2(N)<\ldots<i_{|\Rr(N)|}(N).
\]
Similarly, we denote $\Cc(N)=\{j_1(N),\dots,j_{|\Cc(N)|}(N)\}$ with
\[
j_1(N)<j_2(N)<\ldots<j_{|\Cc(N)|}(N).
\]
Given two sub-matrices $N_1$ and $N_2$, we say that $N_1\preceq_\Rr N_2$ if the rows of $N_1$ come ``before'' the rows of $N_2$; more rigorously,
\begin{multline*}
N_1\preceq_\Rr N_2\ \stackrel{\text{def}}{\Longleftrightarrow}\ |\Rr(N_1)|\le|\Rr(N_2)|\ \text{ and }\ i_k(N_1)\le i_k(N_2)\quad\forall\ 1\le k\le|\Rr(N_1)|.
\end{multline*}
We also give the analogous definition for columns:
\begin{multline*}
N_1\preceq_\Cc N_2\ \stackrel{\text{def}}{\Longleftrightarrow}\ |\Cc(N_1)|\le|\Cc(N_2)|\ \text{ and }\ i_k(N_1)\le i_k(N_2)\quad\forall\ 1\le k\le|\Cc(N_1)|.
\end{multline*}
Note that $\preceq_\Rr$ and $\preceq_\Cc$ are preorders (i.e.\ reflexive and transitive) but not partial orders in general. We write $N_1\prec_\Rr N_2$ if $N_1\preceq_\Rr N_2$ and $N_1\neq N_2$. We define $N_1\prec_\Cc N_2$ similarly.
\bigbreak

It turns out that, when restricted to the set of square submatrices of rank equal to rank $M$, the relations $\preceq_\Rr$ and $\preceq_\Cc$ admit a unique common minimal element, in a sense made precise by the following lemma.

\begin{lem}
\label{lemma:Mstar}
Let $M$ be a $m\times n$-matrix. There exists a submatrix $M^*$ of $M$ such that $\rank M^*=\rank M$ and which is minimal in the following sense: any submatrix $N$ with $\rank N=\rank M$ satisfies
\begin{equation}
\label{eq:minimal}
M^*\preceq_\Rr N\quad\text{ and }\quad M^*\preceq_\Cc N.
\end{equation}
The submatrix $M^*$ is uniquely defined by this condition, and is a square matrix.
\end{lem}

It is easy to check that if such a submatrix $M^*$ exists, then it is unique and square. Thus we just focus on the existence. We will prove existence in a somehow constructive way, by giving a procedure for finding $M^*$. In fact, we will give two different procedures and we will show that they yield the same submatrix; as a consequence, this submatrix will satisfy the conditions demanded to $M^*$.

As a first intermediate step, let us describe two constructions which allow to associate to $M$ two special (non-square in general) submatrices. We call these two sub-matrices $V(M)$ and $H(M)$. Here $V$ stands for \emph{vertical} and $H$ for $\emph{horizontal}$.

Let us first describe how to construct $V(M)$: it is uniquely defined by the properties
\begin{align*}
\Cc(V(M))=&\ \Cc(M)=\{1,\dots,n\}
\\
i\in \Rr(V(M))\Leftrightarrow& \text{ the $i$-th row of $M$ is linearly independent from}
\\
&\text{  the first $i-1$ rows of $M$.}
\end{align*}
(If $i=0$, we mean that $i\in\Rr(V(M))$ if and only if the first row is not identically zero.) This procedure ensures that $V(M)$ is minimal with respect to rows among submatrices of $M$ of maximal rank. More precisely, $V(M)$ satisfies
\[
\rank V(M)=\rank M\quad\text{and}\quad V(M)\preceq_\Rr N\quad\forall\ N\text{ submatrix of }M\text{ with }\rank N=\rank M.
\]
The construction of $H(M)$ is the same as for $V(M)$, but with the roles of rows and columns inverted. More precisely,
\begin{align*}
\Rr(H(M))=&\ \Rr(M)=\{1,\dots,m\}
\\
j\in \Cc(H(M))\Leftrightarrow& \text{ the $j$-th row of $M$ is linearly independent from }
\\
&\text{ the first $j-1$ rows of $M$.}
\end{align*}
Analogously,
\[
\rank H(M)=\rank M\quad\text{and}\quad H(M)\preceq_\Cc N\quad\forall\ N\text{ submatrix of }M\text{ with }\rank N=\rank M.
\]
Now that we have introduced the two constructions, we can iterate them. 
In particular, we can consider  $HV(M)$ and $VH(M)$. We regard them as submatrices of $M$. By construction, they are square non-singular submatrices of size equal to $\rank M$. In fact, it turns out that they coincides, and they are the submatrix $M^*$ which we are looking for. More precisely, the following two claims are true:
\begin{itemize}
\item[(i)] $HV(M)=VH(M)$;
\item[(ii)] the matrix $M^* := HV(M)=VH(M)$ satisfies the conditions required in the statement.
\end{itemize}

Proof of (i). 
We want to prove that $\Rr(HV(M))=\Rr(VH(M))$ and $\Cc(HV(M))=\Cc(VH(M))$. As already pointed out, $HV(M)$ and $VH(M)$ are square submatrices of equal size (equal to $\rank M$). Hence it suffices to prove that $\Rr(VH(M))\subseteq \Rr(HV(M))$ and that $\Cc(HV(M))\subseteq \Cc(VH(M))$. We focus just on the first inclusion, the second being analogous. 

By the properties of the constructions $H$ and $V$ described above, we have:
\begin{align*}
k\in\Rr(VH(M)) \Rightarrow &\text{ the $k$-th row of $H(M)$ is linearly independent }
\\
& \text{ from the first $k-1$ rows of $H(M)$}
\\
\Rightarrow & \text{ the $k$-th row of $M$ is linearly independent }
\\
& \text{ from the first $k-1$ rows of $M$}
\\
\Rightarrow &\ k\in \Rr(V(M))=\Rr(HV(M))
\end{align*}
as desired.

Proof of (ii). Let $N$ be a sub-matrix of $M$ with $\rank N=\rank M$. By the properties of $V(M)$, $V(M)\preceq_{\Rr} N$. From $M^*=HV(M)$ we deduce $\Rr(M^*)=\Rr(HV(M))=\Rr(V(M))$ and thus $M^*\preceq_{\Rr}N$ as well. The proof of $M^*\preceq_{\Cc}N$ is similar.
\renewcommand\qedsymbol{\ensuremath{\openbox}~\textsf{\textsc{\footnotesize{Lemma \ref{lemma:Mstar}}}}~\ensuremath{\openbox}}\qedhere

Let us emphasize the following characterization of $M^*$ which follows directly from the proof above:
\begin{equation}
\label{eq:property Mstar}
\begin{aligned}
i\in \Rr(M^*)\Leftrightarrow& \text{ the $i$-th row of $M$ does not belong to the linear span }
\\
&\text{ of the first $i-1$ rows of $M$.}
\\
j\in \Cc(M^*)\Leftrightarrow& \text{ the $j$-th column of $M$ does not belong to the linear span }
\\
&\text{ of the first $j-1$ columns of $M$.}\renewcommand\qedsymbol{}\qedhere
\end{aligned}
\end{equation}

\textbf{End of the proof of Proposition \ref{prop:core proposition}.}
We fix $r$ between $m$ and $\dim J^p(X,Y)-1$ 
 and assume that $Z^r_{a,V}$ is not empty (otherwise we have nothing to prove).
 Let $z_0$ be an arbitrary element of $Z^r_{a,V}$. Our goal will now be to find $\theta$ different square submatrices
$[M_{a,V}(z)]_i$ of size $r+1$ in $M_{a,V}(z)$ such that the equations 
$\det [M_{a,V}(z)]_i=0$ are independent near $z_0$.
The functions $z\lmto \det [M_{a,V}(z)]_i$ are then the components of the  map 
$F_{a,V}(z)$ mentioned in Proposition \ref{prop:core proposition}.
These functions are clearly polynomials in $z$, with coefficients depending
smoothly on $a$ and $V$.
We omit from now on to explicitly mention $a$ and $V$, and note $M(z)$ instead of 
$M_{a,V}(z)$.

Let $M^*_{z_0}$ be the square submatrix of size $r$ associated to $M(z_0)$ by the Lemma \ref{lemma:Mstar}. 
We take the notational convention that the symbol $M^*_{z_0}$ 
without further specifications stands for a pattern of rows and columns, i.e.\ $M^*_{z_0}$ is the datum $(\Rr(M^*_{z_0}),\Cc(M^*_{z_0}))$. We can also identify
 $M^*_{z_0}$ to a matrix-valued function of $z$, but in this case we explicitly write 
$M^*_{z_0}(z)$ or $M^*_{z_0}(\cdot)$. This is in order to avoid ambiguities and distinguish, 
for instance, between $M^*_{z_0}$ and $M^*_{z_0}(z_0)$.
 We adopt the same convention for all the submatrices encountered below,
 such as $\hat M_{z_0}, M^*_{z_0,\#(i,j)}$, etc., which we shall define shortly.

We have 
\[
\deter M^*_{z_0}(z_0)\neq 0\quad\text{and}\quad\rank M^*_{z_0}(z_0)=\rank M(z_0)=r.
\]
Let us call $\hat M_{z_0}$ the submatrix whose rows and columns are exactly the ones \emph{not} appearing in $M^*_{z_0}$, i.e.\ 
\[
\Rr(\hat M_{z_0})=\{1,\dots,\dim J^p(X,Y)\}\setminus \Rr(M^*_{z_0}),\qquad \Cc(\hat M_{z_0})=\{1,\dots,m+n\}\setminus \Cc(M^*_{z_0})
\]
As already mentioned above, the first $m$ columns of $M$ are linearly independent as well as the first $\dim J^{p-1}(X,Y)$-rows. By the characterization \eqref{eq:property Mstar} we deduce that $\hat M_{z_0}$ is entirely contained in the bottom-right block of $M$, i.e.\ it is a submatrix of $B$:
\begin{align*}
\Rr(\hat M_{z_0})\subseteq \Rr(B)&=\{\dim J^{p-1}(X,Y)+1,\dots,\dim J^p(X,Y)\}
\\
\Cc(\hat M_{z_0})\subseteq \Cc(B)&=\{m+1,\dots,m+n\}.
\end{align*}
Let us denote $M^*_{z_0,\#(i,j)}$ the submatrix obtained by adding to $M^*_{z_0}$ the $i$-th row and the $j$-th column of $M$, i.e.
\[
\Rr(M^*_{z_0,\#(i,j)})=\Rr(M^*_{z_0})\cup\{i\},\qquad\Cc(M^*_{z_0,\#(i,j)})=\Cc(M^*_{z_0})\cup\{j\}.
\]
We are interested to the case when $(i,j)$ belongs to $\Rr(\hat M_{z_0})\times\Cc(\hat M_{z_0})$. In this case the submatrix $M^*_{z_0,\#(i,j)}$ is a square submatrix of size $r+1$. We are in the following situation:
\[
\begin{cases}
\deter M^*_{z_0}(z_0)\neq 0
\\
\deter M^*_{z_0,\#(i,j)}(z_0)= 0\qquad \forall \ (i,j)\in\Rr(\hat M_{z_0})\times\Cc(\hat M_{z_0})
\end{cases}
\]
Proposition \ref{prop:core proposition} follows from the following two lemma:

\begin{lem}
\label{triangular}
If $(i_1,j_1),\dots,(i_\theta,j_\theta)$ are pairwise distinct couples in $\Rr(\hat M_{z_0})\times\Cc(\hat M_{z_0})$ such that
\begin{equation}
\label{eq:increasing indexes}
i_1\le i_2\le \dots\le i_\theta\quad\text{and}\quad j_1\le j_2\le \dots\le j_\theta,
\end{equation}
then the differentials evaluated at $z_0$
\[
d_{z_0}\deter M^*_{z_0,\#(i_1,j_1)}(\cdot),\,\dots,\,d_{z_0}\deter M^*_{z_0,\#(i_\theta,j_\theta)}(\cdot)
\]
are linearly independent.
\end{lem}

\begin{lem}\label{theta}
 For $\theta=n+1-c +2 \big(\dim J^p(X,Y)-1-r\big)$, 
there do exist pairwise distinct couples as above.
\end{lem}

\textsc{Proof of Lemma \ref{theta}.}
One may for instance consider 
the couples of indexes successively encountered along the following ``path'' 
in the matrix $\hat M_{z_0}$: starting from the upper-left corner of the matrix, 
and then moving horizontally along the first row until the upper-right corner,
 and then moving vertically along the last column until the bottom-right corner. 
It is clear that the couples of indexes successively encountered along this path 
satisfy the condition \eqref{eq:increasing indexes}. Their number is the ``semi-perimeter'' 
of the matrix, or more rigorously
$
\big|\Rr(\hat M_{z_0})\big|+\big|\Cc(\hat M_{z_0})\big|-1.
$
Recalling the definition of $\hat M_{z_0}$, this is the same as
\begin{align*}
|\Rr(M)|-|\Rr(M^*_{z_0})|+|\Cc(M)|-|\Cc(M^*_{z_0})|-1&=\left(\dim J^p(X,Y)-r\right)+\left(m+n-r\right)-1
\\
&=n+1-c +2 \big(\dim J^p(X,Y)-1-r\big)
\end{align*}
This ends the proof of Lemma \ref{theta}.
\qed

\textsc{Proof of Lemma \ref{triangular}.}
 For $(i,j)\in\Rr(B)\times\Cc(B)$, 
we recall that $i$ is the index of a line of $B$, hence it corresponds to a coordinate
$y^{s(i)}_{\alpha(i)}$ of $\F_{p-1}^p$, while $j$ is an integer 
%between $1$ and $n$
between $m+1$ and $m+n$. 
Then, the coefficient of the matrix $B(z)$ at line $i$ and column $j$ is just
$y^{s(i)}_{\alpha(i)+\delta_{j-m}}$. It is a component, that we denote by $z_{[i,j]}$, of $z$.
Note however that the same component of $z$ 
may appear at several different places in the matrix $B(z)$. 
It can happen that $z_{[i,j]}=z_{[i',j']}$ with $(i,j)\neq (i',j')$,
it is the case  when $s(i)=s(i')$ and $\alpha(i)+\delta_{j-m}=\alpha(i')+\delta _{j'-m}$.
Our order on the coordinates allows us to overcome this difficulty
thanks to the following Lemma:

\begin{lem}
\label{lexicographic property} 
Let $(i,j)$ and $(h,k)$ belong to $\Rr(B)\times\Cc(B)$, and satisfy
\[
h\ge i,\ k\ge j,\ (i,j)\neq (h,k).
\]
Then, $z_{[h,k]}\neq z_{[i,j]}$ and, if $(i,j)$ and $(h,k)$ belong to $\Rr(\hat M_{z_0})\times\Cc(\hat M_{z_0})$, then
\begin{align*}
\frac{\partial}{\partial z_{[i,j]}}\deter M^*_{z_0,\#(i,j)}(z_0)&=\pm\deter M^*_{z_0}(z_0)\neq0
\\
\frac{\partial}{\partial z_{[h,k]}}\deter M^*_{z_0,\#(i,j)}(z_0)&=0.
\end{align*}
\end{lem}

Lemma  \ref{lexicographic property} implies that 
$
z_{[i_1,j_1]},\dots,z_{[i_\theta,j_\theta]}
$
are pairwise distinct components of $z$, and that the square matrix\[
%\displaystyle
\left[
\begin{matrix}
\displaystyle\frac{\partial}{\partial z_{[i_1,j_1]}}\deter M^*_{z_0,\#(i_1,j_1)}(z_0)  
&\qquad  \dots\qquad & \displaystyle\frac{\partial}{\partial z_{[i_\theta,j_\theta]}}\deter M^*_{z_0,\#(i_1,j_1)}(z_0)
\cr
\vdots & \ddots & \vdots
\cr
\displaystyle\frac{\partial}{\partial z_{[i_1,j_1]}}\deter M^*_{z_0,\#(i_\theta,j_\theta)}(z_0)
  &\qquad  \dots\qquad & \displaystyle\frac{\partial}{\partial z_{[i_\theta,j_\theta]}}\deter M^*_{z_0,\#(i_\theta,j_\theta)}(z_0)
\cr
\end{matrix}
\right]
\]
has the form
\begin{equation}
\label{eq:triangular}
\left[\begin{matrix}
\pm\deter M^*_{z_0}(z_0) & 0 & 0 & \ldots & 0
\cr
* & \pm\deter M^*_{z_0}(z_0) & 0 & \ldots & 0
\cr
* & * & \pm\deter M^*_{z_0}(z_0) & \ldots & 0
\cr
\vdots & \vdots & \vdots & \ddots & \vdots
\cr
* & * & * & \ldots & \pm\deter M^*_{z_0}(z_0)
\cr
\end{matrix}
\right]
\end{equation}
hence is  invertible, since $\deter M^*_{z_0}(z_0)\neq 0$.
 This proves Lemma \ref{triangular}, using Lemma \ref{lexicographic property}.
\qed

\textsc{Proof of Lemma \ref{lexicographic property}.}
Let us first prove that $z_{[i,j]}\neq z_{[h,k]}$.
If $s(i)\neq s(h)$, then the conclusion holds.
If $s(i)= s(h)$, then $\alpha(i)\geq \alpha(h)$
for the lexicographic order.
On the other hand, the inequality $k\geq j$ implies that
%$\delta_{k}\leq \delta_{j}$ 
$\delta_{k-m}\leq \delta_{j-m}$ 
 for the lexicographic order.
These two inequalities do not sum to an equality
because they are not both equalities
(recall the hypothesis $(i,j)\neq (h,k)$), hence
%$\alpha(i)+\delta_{j}\neq\alpha(h)+\delta_{k}$,
$\alpha(i)+\delta_{j-m}\neq\alpha(h)+\delta_{k-m}$,
 and then
$z_{[i,j]}\neq z_{[h,k]}$.

To prove the equality
\begin{equation}
\label{eq:partial derivative minor}
\frac{\partial}{\partial z_{[h,k]}}\deter M^*_{z_0,\#(i,j)}(z_0)=0,
\end{equation}
let us consider, for every $(h',k')\in \Big(\Rr(M^*_{z_0,\#(i,j)})\times\Cc(M^*_{z_0,\#(i,j)})\Big)$, the submatrix $N_{h',k'}$ defined by
\begin{align*}
\Rr(N_{h',k'})=\Rr(M^*_{z_0,\#(i,j)})\setminus\{h'\}=\Big(\Rr(M^*_{z_0})\cup\{i\}\Big)\setminus\{h'\}
\\
\Cc(N_{h',k'})=\Cc(M^*_{z_0,\#(i,j)})\setminus\{k'\}=\Big(\Cc(M^*_{z_0})\cup\{j\}\Big)\setminus\{k'\}.
\end{align*}
We have
\[
\frac{\partial}{\partial z_{[h,k]}}\deter M^*_{z_0,\#(i,j)}(z_0)=\sum_{(h',k')}\pm\deter N_{h',k'}(z_0)
\]
where the actual sign  is irrelevant  and the sum is taken over all couples $(h',k')\in \Big(\Rr(M^*_{z_0,\#(i,j)})\times\Cc(M^*_{z_0,\#(i,j)})\Big)$  such that $z_{[h',k']}=z_{[h,k]}$.

We claim that each square matrix $N_{h',k'}$ is singular, thus proving \ref{eq:partial derivative minor}. In view of the definition of $M^*_{z_0}$ in 
Lemma \ref{lemma:Mstar}, it is enough to observe that we can't have 
both $M^*_{z_0}\preceq_\Rr N_{h',k'}$ and $ M^*_{z_0}\preceq_\Cc N_{h',k'}$.
This would imply that we have both 
$h'\leq i$ and $k'\leq j$ and then that 
 $h'\leq h$ and $k'\leq k$.
As we have already seen, since  $z_{[h',k']}=z_{[h,k]}$, this would imply that $(h',k')=(h,k)$. Since $h'\le i\le h$ and $k'\le j\le k$, we would finally have $(h,k)=(i,j)$, in contradiction with our hypotheses.

Finally, we have
$$
\frac{\partial}{\partial z_{[i,j]}}\deter M^*_{z_0,\#(i,j)}(z_0)=\pm \deter M^*_{z_0}(z_0)+\sum_{(i',j')}\pm\deter N_{i',j'}(z_0)
$$
where the sum is taken on 
all couples $(i',j')\in \Big(\Rr(M^*_{z_0,\#(i,j)})\times\Cc(M^*_{z_0,\#(i,j)})\Big)$ 
such that $z_{[i',j']}=z_{[i,j]}$.
We conclude as above that all the terms in the sum vanish.
\qed

\subsection{The analytic case}\label{Aproof}
We prove Theorem \ref{AT} by recurrence on $p$.
When $p=0$, $A=A_1$ hence the statement follows from Theorem \ref{NDT}.

Since $A_0$ is defined by analytic conditions (at least in a suitable chart), it is a stratified set. It suffices to bound the codimension of the stratum $S\subset A_0$ of maximal dimension. Let us consider the restricted projection 
\[
\big({\pi^p_{p-1}}_{|S}\big)\colon S\to J^{p-1}(X,Y).
\]
and the associated rank map
\[
S\ni a\mapsto \rank d_a\,\big({\pi^p_{p-1}}_{|S}\big)\quad\in\ \Big\{0,1,\dots,\min\{\dim S,\dim J^{p-1}(X,Y)\}\Big\}.
\]
Let us also consider an open subset $U$ of $S$ such that the rank map is constant on $U$. Such a subset exists, for instance we can take as $U$ the preimage of the maximum value attained by the map (this preimage is open because the rank map is lower-semicontinuous).

It follows from the constant-rank theorem  that, up to further restricting $U$ if necessary, $\pi^p_{p-1}(U)$ is a submanifold of $J^{p-1}(X,Y)$. Let us call $V$ this manifold. We claim that
\[
U\subseteq \tilde V
\]
where $\tilde V\subseteq J^p(X,Y)$ is defined according to \eqref{eq:tilde}, i.e.\ 
\[
\tilde V=\big\{ j^{p}_xf\in J^{p}(X,Y): j^{p-1}f \text{ is not transverse to }V\text{ at }x \big\}.
\]
Since the conclusion of Theorem \ref{AT} is assumed to be true for $p-1$, we have $\codim \tilde V\ge n+1$, hence the claim implies $\codim U\ge n+1$. Since $U$ is open in $S$ and $S$ is the stratum of maximal dimension, we get
\[
\codim A_0=\codim S=\codim U\ge n+1,
\]
which proves the proposition.
 Let us now prove the claim $U\subseteq\tilde V$. Given any $a=j^p_xf\in U$, we have 
\[
T_{\pi^p_{p-1} a}V=\big({\pi^p_{p-1}}_{|S}\big)(T_a U)\subseteq \pi^p_{p-1}(T_a A).
\]
Here the first equality follows by the constant-rank theorem, while the inclusion follows from the fact that $U\subseteq S\subseteq A$. Moreover, the very definition of $A_0$ yields
\[
\pi^p_{p-1}(T_a A)+E^{p-1}(a)\subsetneq T_{\pi^p_{p-1} a} J^{p-1}(X,Y).
\]
It follows that
\[
T_{\pi^p_{p-1} a}V+E^{p-1}(a)\subseteq \pi^p_{p-1}(T_a A)+E^{p-1}(a) \subsetneq T_{\pi^p_{p-1} a} J^{p-1}(X,Y)
\]
which implies that $a\in \tilde V$, as desired.
\qed

\end{document}